\input ./style/arxiv-general.cfg
\documentclass[aap,MSNbibl,dvips]{arximspdf}
\makeatletter
   \@ifpackageloaded{graphicx}{}{\usepackage{graphicx}}
\makeatother

\doi{10.1214/15-AAP1127}
\volume{26}
\issue{3}
\pubyear{2016}
\firstpage{1620}
\lastpage{1635}
\docsubty{FLA}

\makeatletter
\let\epsilon\varepsilon
\newcommand{\eqref}[1]{(\ref{#1})}
\newtheorem{lemma}[thm]{Lemma}
\newtheorem{prop}[thm]{Proposition}

\newproclaim{define}[thm]{Definition}
\newproclaim{conjecture}[thm]{Conjecture}
\newproclaim{question}[thm]{Open question}

\newproclaim{fact}{Fact}
\newproclaim{remark}[thm]{Remark}
\makeatother

\begin{document}
\begin{frontmatter}

\title{From transience to recurrence with Poisson tree frogs}
\runtitle{From transience to recurrence with Poisson tree frogs}

\begin{aug}
\author[A]{\fnms{Christopher}~\snm{Hoffman}\thanksref{M1,T1}\ead
[label=e1]{hoffman@math.washington.edu}},
\author[B]{\fnms{Tobias}~\snm{Johnson}\thanksref{M2,T2}\corref{}\ead[label=e2]{tobias.johnson@usc.edu}}\\
\and
\author[A]{\fnms{Matthew}~\snm{Junge}\thanksref{M1}\ead[label=e3]{jungem@math.washington.edu}}
\runauthor{C. Hoffman, T. Johnson and M. Junge}
\affiliation{University of Washington\thanksmark{M1} and University of
Southern California\thanksmark{M2}}
\address[A]{C. Hoffman\\
M. Junge\\
Department of Mathematics\\
University of Washington\\
Seattle, Washington 98195\\
USA\\
\printead{e1}\\
\phantom{E-mail:\ }\printead*{e3}}
\address[B]{T. Johnson\\
Department of Mathematics\\
University of Southern California\\
Los Angeles, California 90089\\
USA\\
\printead{e2}}
\end{aug}
\thankstext{T1}{Supported in part by NSF Grant DMS-13-08645 and NSA
Grant H98230-13-1-0827.}
\thankstext{T2}{Supported in part by NSF Grant DMS-14-01479.}

%
\received{\smonth{1} \syear{2015}}
%
\revised{\smonth{6} \syear{2015}}

%
\begin{abstract}
Consider the following interacting particle system on the $d$-ary tree,
known as the \emph{frog model}:
Initially, one particle is awake at the root and i.i.d. Poisson many
particles are sleeping at every
other vertex.
Particles that are awake perform simple random walks, awakening
any sleeping particles they encounter.
We prove that there is a phase transition between transience and
recurrence as
the initial density of particles increases, and we give the order of
the transition up to a logarithmic factor.
\end{abstract}

%
\begin{keyword}[class=AMS]
\kwd{60K35}
\kwd{60J80}
\kwd{60J10}
\end{keyword}
\begin{keyword}
\kwd{Frog model}
\kwd{transience}
\kwd{recurrence}
\kwd{phase transition}
\end{keyword}
\end{frontmatter}

\section{Introduction}


We study a system of branching random walks
known as the frog model, and we discover a phase transition as the
initial state
becomes more saturated with particles.
Similar phase transitions have been observed in related models, including
activated random walk \cite{DRS,ST}, reinforced random walk \cite
{reinforcedrw}, 
killed branching random walk \cite{killedbrw} and the contact process
\cite{contact}.

The frog model starts
with a single particle awake at the root of a graph and sleeping
particles at the other vertices.
The initial configuration of sleeping particles can be deterministic or
random. Particles that are awake perform independent simple random
walks in discrete time.
When a vertex with sleeping particles is first visited, all of the
particles at the site wake up and each begins its own walk. The name
``frog model'' was coined in 1996 by Rick Durrett; we continue the
zoomorphism and refer to the particles as frogs.
As with other interacting particle systems, the frog model is often
motivated as a model
for the spread of a rumor or infection (see \cite{shape}, e.g.). It and its
variants have also
found interest as models of combustion \cite{combustion,CQR,RS},
generally with particles moving in continuous time.

We call a realization of the frog model \emph{recurrent} if the root is
visited infinitely
often by frogs and \emph{transient} if not.
Even if each individual frog is transient, the aggregate of visits to
the root can still be infinite.
For this reason, the transience or recurrence of the frog model gives a
measurement of its growth, and
the question of transience or recurrence for the frog model on a given
graph is one
of the most fundamental ones.

The first ever published result on the frog model is that it is
recurrent on $\mathbb{Z}^d$ with one sleeping
frog per site for all $d$ \cite{telcs1999}. In fact, the frog model on
$\mathbb Z^d$ is recurrent for any i.i.d. initial
configuration of sleeping frogs \cite{randomshape}.
It is natural to wonder if a sparser configuration changes the behavior.
\cite{randomrecurrence} exhibits a threshold at which a frog model
with $\operatorname{Bernoulli}(\alpha\|x\|^{-2})$ frogs at each $x \in\mathbb{Z}^d$
switches from transience to recurrence.
A similar phenomenon occurs
when the walks have a bias in one direction: \cite{recurrence} finds
that on $\mathbb{Z}$,
the model is recurrent if and only if the number of sleeping frogs per
site has infinite logarithmic moment.
Recently, this result was partially extended to $\mathbb{Z}^d$ in \cite
{newdrift} and worked out in finer detail in \cite{integersdrift2}.

Let $\mathbb{T}_d$ denote the full infinite $d$-ary tree, in which the
root has degree~$d$ and all
other vertices degree~$d+1$.
The question of transience or recurrence on
$\mathbb{T}_d$ is especially subtle. On one hand, the number of sleeping frogs
grows exponentially with the distance from the root.
On the other hand, each frog that wakes up has a drift away from the root;
its probability of visiting the root shrinks exponentially as the
starting vertex of the frog
moves outward.
The question of whether $\mathbb{T}_d$ is transient for the one-per-site
model is posed in \cite{phasetree} and again in \cite{frogs} and
\cite
{recurrence}. Surprisingly, the answer depends on the degree of the
tree. In \cite{HJJ}, we prove that the one-per-site frog model is
recurrent on the binary tree and transient on $d$-ary trees with $d
\geq5$.

We conjecture that the one-per-site frog model is recurrent for $d=3$
and transient for $d=4$.
While we would like to pin this down and complete the picture of
transience and recurrence
for the one-per-site frog model on trees, we believe that the most
interesting aspect of this
work is that the frog model on trees is teetering on the edge between
recurrence and
transience. The point of this paper is to demonstrate this more precisely.
We consider the frog model on $\mathbb{T}_d$ with i.i.d. $\operatorname
{Poi}(\mu)$
sleeping frogs at each site.
Our result is a phase transition between recurrence and transience as
$\mu$ varies.


\begin{thm}\label{thm1}
Consider the frog model on a $d$-ary tree with $\operatorname{Poi}(\mu)$
sleeping frogs per site. For all $d\geq2$, there exists a critical
value $\mu_c(d)>0$ such that the model is recurrent a.s. if $\mu>\mu
_c(d)$ and transient a.s. if $\mu<\mu_c(d)$. The critical value satisfies
\[
Cd < \mu_c(d) < C'd \log d
\]
%
for some constants $C$ and $C'$.
\end{thm}

\begin{pf}
By a straightforward coupling, the probability of recurrence is
monotone in $\mu$.
By \cite{HJJ}, Theorem~4, the probability of recurrence is either $0$
or $1$.
The theorem is then an immediate consequence of Propositions~\ref
{proprecurrence} and~\ref{proptransience}, where
we prove recurrence and transience, respectively.
\end{pf}

Contrast our result with the frog model on $\mathbb{Z}^d$, which is
recurrent
for any i.i.d. configuration of sleeping frogs \cite{randomshape}.
To show the existence of the recurrence phase, we consider a restricted
process that lets us take advantage
of the recursive structure of $\mathbb{T}_d$. We then use a bootstrapping
argument, showing that the number
of returns to the root is stochastically larger and larger at each step.
We establish the transience phase essentially by dominating the model
with a branching random walk, using
a similar argument as in \cite{HJJ}.
As in that paper, the most difficult part is recurrence.
Our result is an advance in that we are able to show recurrence
on any $d$-ary tree with enough sleeping frogs.
In \cite{HJJ}, we prove recurrence only for $d=2$, and the
proof does not apply to a general choice of $d$; even
extending it to $d=3$ seems difficult. The argument here relies on
having Poisson many sleeping frogs at each site, however, and thus
neither result implies the other. A more detailed comparison between
the recurrence proofs in the two papers is in Section~\ref{subseccomparison}.

\subsection*{Further questions}

A nice general survey on the frog model can be found in \cite{frogs}.
Here, we pose four questions specifically related to the frog model on trees.

The question most directly related to our paper is to better estimate
the critical value $\mu_c(d)$. We are interested in both the asymptotic
behavior and precise values for small $d$.

\begin{question}
What is the correct order of $\mu_c(d)$ as $d \to\infty$? Also, what
is the value of $\mu_c(d)$ for small $d$?
\end{question}

We suspect that $\mu_c(d)=\Theta(d)$.
As for the second question,
the best bounds we can prove for $d=2$ are
$0.125 \leq\mu_c(2) \leq1.13$ (see Section~\ref{subseccomparison}).

As a start at considering the frog model on less regular graphs, we
would like to know if
the analogue of our result holds on Galton--Watson trees.

\begin{question}
Consider a frog model with $\operatorname{Poi}(\mu)$ frogs at each
site of an
infinite Galton--Watson tree. As $\mu$ varies, does a phase transition
occur between transience and recurrence?
\end{question}

We are also interested in the relationship between the frog
model and the degree distribution of the tree.

\begin{question}
Does the recurrence of the frog model on a Galton--Watson tree depend
on the entire degree distribution
or just the maximal degree?
Concretely, consider a one-per-site frog model on a Galton--Watson
tree where each vertex
has probability $p$ of having two children and probability $1-p$ of
having five children.
\cite{HJJ}, Theorem~1 implies that this is recurrent when $p=1$ and
transient when $p=0$.
Is it recurrent for any $p< 1$?
\end{question}

This dependence on the maximal degree of the tree alone
is seen in the contact process 
(see \cite{contact} and \cite{pemantle2001branching}, Proposition~2.5).

Our next question comes from Itai Benjamini and concerns the frog model
on finite trees. Define the \emph{cover time} to be the expected time
for every frog to wake up in a one-per-site frog model on the full
$d$-ary tree with height~$n$. We call this the cover time since it is
equivalent to the time when every site is visited.
A naive bound on the cover time is $O(n^2 d^n)$, the expected time for
a single random walk to visit every site, as shown in \cite{aldous}. We
have an unpublished proof improving this to $O(n^5( {d}/{\sqrt2}
)^n)$, but we suspect the correct value is polynomial.

\begin{question}
Is the cover time for the one-per-site frog model on a $d$-ary tree
of height~$n$ polynomial in $n$?
\end{question}

Possibly the cover time on finite trees relates to the recurrence and
transience properties on the corresponding infinite tree. For instance,
it would be exciting to see that the cover time
is polynomial in the height of the tree for $d=2$ but exponential for
higher $d$.
This would be reminiscent of the contact process, which behaves
similarly on finite lattices and
trees as on their infinite counterparts
\cite{contactfinite1,contactfinite2,contactfinite3,CMMV}.

\section{Recurrence}
\label{secrecurrence}


We start with a sketch. Let $\nu'$ be the law of the number of visits
to the root in the frog model with $\operatorname{Poi}(\mu)$ frogs at
each site.
To get some regularity, we restrict the motion of awakened frogs to
the nonbacktracking component of their ranges. Call this the \emph
{nonbacktracking frog model} (more details are in Section~\ref{secnonbacktracking}) and let $\nu$ be the law of the number of visits
to the root in this model. A coupling argument in Proposition~\ref
{propbacktrackingcoupling} confirms the intuition that
\begin{eqnarray}
\label{eqdominance} &&\nu\preceq\nu'.
\end{eqnarray}
Here, $\preceq$ denotes \emph{stochastic dominance}, that is
$ \nu([x, \infty)) \leq\nu'( [ x ,\infty))$ for all $x$.

%
\begin{figure}

\includegraphics{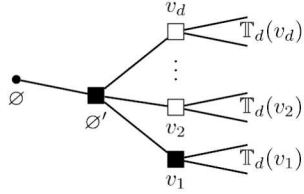}

\caption{The frog from $\varnothing$ visits $\varnothing'$ and
$v_1$. Suppose %
at most one frog in the nonbacktracking frog model is allowed to enter
each $\mathbb{T}_d(v_i)$ and
only frogs woken at $\varnothing'$ and emerging from $\mathbb
{T}_d(v_1)$ can
enter other subtrees. We see in Lemma~\protect\ref{lemfixedpoint}
that the number
of visits to $\varnothing$ is stochastically fewer than $\nu$ and is
distributed as $\mathcal{A}\nu$.}\label{figtree}
\end{figure}

In Section~\ref{secAdef}, we define an operator $\mathcal A$ under
which the image of $\nu$ has an interpretation in an even more
restricted frog model. 
First, a bit of notation (see Figure~\ref{figtree}) is necessary.
Say the initial nonbacktracking frog
moves down the tree from the root $\varnothing$
to $\varnothing'$ and then to $v_1$. Let $v_2,\ldots,v_d$ be the
other children
of $\varnothing'$ and let $\mathbb{T}_d(v_i)$ denote the subtree
rooted at $v_i$.
The measure $\mathcal A \nu$ is the law of the number of visits to the
root in the nonbacktracking frog model with two further restrictions:
\begin{longlist}[(ii)]
\item[(i)] At most one frog can enter $\mathbb{T}_d(v_i)$ for
each $1
\leq i \leq d$.
\item[(ii)] Only frogs woken at $\varnothing'$ and those emerging
from $\mathbb{T}_d(v_1)$ can enter the other $\mathbb{T}_d(v_i)$.
\end{longlist}
The advantage of (i) is that it makes the number of frogs
emerging from the activated subtrees i.i.d.
random variables. The advantage of (ii) is that it simplifies
which subtrees become activated (see Lemma~\ref{lemfirstwavecondition}).
Intuitively, these restrictions reduce the number of visits to the root.
This is made rigorous in Lemma~\ref{lemfixedpoint} where we prove that
\begin{eqnarray}
\label{eqcontraction}
&& \mathcal{A} \nu\preceq\nu.
\end{eqnarray}
We stress that this a special property of $\nu$.
In fact, the essence of our argument is to show that when $\mu$ is
large enough,
\eqref{eqcontraction} can hold only if $\nu=\delta_\infty$.


Section~\ref{secAproperties} explores properties of $\mathcal A$. In
Lemma~\ref{lemmonotonicity}, we show that $\mathcal A$ is monotonic,
meaning that for two probability measures $\pi_1$ and $\pi_2$,
\begin{eqnarray}\label{eqmonotonic}
&&\mbox{if $\pi_1 \preceq\pi_2$, then $\mathcal A
\pi_1 \preceq \mathcal A \pi_2$.}
\end{eqnarray}
Lemma~\ref{lemmixture} shows that $\mathcal A$ acts nicely on the Poisson
distribution. In fact, by writing the Poisson distribution in a
nonstandard way (see Lemma~\ref{lemPoidivision}), we can compare
$\mathcal
A \operatorname{Poi}(\lambda)$ with $\operatorname{Poi}(\lambda+
\epsilon)$. We carry this out in
Proposition~\ref{proppoisson},
where we show that when $\mu\geq2 ( d+1) \log d$, there exists
$\epsilon$ such that
\begin{eqnarray}\label{eqpoisson}
&& \operatorname{Poi}(\lambda+ \epsilon) \preceq\mathcal A \operatorname{Poi}(
\lambda)
\end{eqnarray}
for all $\lambda\geq0$.
This is where the value of $\mu$ plays a role.
Proving \eqref{eqpoisson} reduces to comparing two binomial
distributions with parameters depending on $\mu$.

Now we explain how \eqref{eqdominance}, \eqref{eqcontraction},
\eqref
{eqmonotonic} and \eqref{eqpoisson} imply the recurrence part of
Theorem~\ref{thm1}.

\begin{prop} \label{proprecurrence}
If $\mu>2(d+1)\log d$, then the frog model is recurrent~a.s. on the
$d$-ary tree with
an initial configuration of
$\operatorname{Poi}(\mu)$ sleeping frogs per vertex.
\end{prop}

\begin{pf}
By \eqref{eqdominance}, it suffices to prove that $\nu$ is a point
mass at infinity. From \eqref{eqcontraction}, we have
\begin{eqnarray*}
&&\operatorname{Poi}(0) \preceq\mathcal{A}\nu\preceq \nu.
\end{eqnarray*}
Statement~\eqref{eqmonotonic} implies this relation is preserved under
iterations of $\mathcal{A}$. Moreover, \eqref{eqpoisson} lets us
increase the
Poisson term by $\epsilon$ with each iteration. In symbols, this says
that for all $n \geq1$,
\[
\operatorname{Poi}(\epsilon n) \preceq\mathcal A^n \nu\preceq
\mathcal A^{n-1}\nu \preceq\cdots\preceq\mathcal A \nu \preceq\nu.
\]
Taking $n\to\infty$ implies that $\nu$ is a point mass at infinity,
and so the frog model is recurrent almost surely.
\end{pf}

In the rest of this section, we will carry out this plan and prove statements
\eqref{eqdominance}--\eqref{eqpoisson}.
First, we give some notation.
Recall that $\preceq$ denotes stochastic domination. We also use the
notation $X\preceq Y$ to indicate that the law of $X$ is
stochastically dominated by the law of $Y$.
An equivalent condition to stochastic dominance is that
$\pi_1\preceq\pi_2$ if and only if there exists a coupling $(X,Y)$ with
$X\sim\pi_1$,
$Y\sim\pi_2$, and $X\leq Y$ a.s. A thorough reference on stochastic
domination is \cite{SS}.

For a nonnegative random variable $N$, we use $\operatorname{Poi}(N)$
to denote a
mixture of Poisson
distributions; when we write $X\sim\operatorname{Poi}(N)$, we mean
that $X$ is
coupled with $N$ such that
the distribution of $X$ conditional on $N=n$ is $\operatorname
{Poi}(n)$. If $N\sim\pi
$, we also
use $\operatorname{Poi}(\pi)$ to denote the same Poisson mixture. We
similarly use
the notation
$\operatorname{Bin}(N, p)$ and $\operatorname{Bin}(\pi,p)$.

\subsection{The nonbacktracking frog model} \label{secnonbacktracking}

A \emph{random nonbacktracking walk} on $\mathbb{T}_d$ starting at a vertex
$x_0$ moves
in its first step to a uniformly random neighbor of $x_0$. In all
subsequent steps,
it moves to a vertex chosen uniformly from all its neighbors except for
the one it just
arrived from.

Suppose that $(S_n, n\geq0)$ is a random nonbacktracking walk
starting from $x_0$, stopped
if it arrives at the root at step~1 or beyond. (If $x_0$ is the root, then
it is never stopped.)
Define the \emph{nonbacktracking frog model} just as the usual frog
model, except that
the motion of a frog waking at $x_0$
is an independent copy of $(S_n)$, rather than a simple random walk.
The advantage is that when a nonbacktracking frog moves away from the
root, it will forever remain in the just-entered subtree. This gives
the model more self-similarity.
As shown in \cite{HJJ}, Proposition~7, $(S_n)$ can be coupled with
a simple random walk on $\mathbb{T}_d$ starting from $x_0$ so that its
path is
a subset of the simple random
walk's path.
This lets us relate the nonbacktracking and usual frog models, proving
\eqref{eqdominance}.

\begin{prop}\label{propbacktrackingcoupling}
Let $\nu$ and $\nu'$ be the laws of the number of returns to the root
in the nonbacktracking and usual frog models on $\mathbb{T}_d$, respectively,
both with $\operatorname{Poi}(\mu)$ sleeping frogs per vertex.
Then $\nu\preceq\nu'$.
\end{prop}
\begin{pf}
It suffices to show that we can couple the two models so that at least
as many
frogs visit the root in the usual model as in the nonbacktracking model.
We construct the coupling as follows. For each vertex $v\in\mathbb{T}_d$,
make the number of sleeping
frogs on $v$ identical in the two models. Make each frog's path in the
nonbacktracking
model a subset of the corresponding frog's path in the usual model as
previously described.
Thus, any frog woken in the nonbacktracking model is also woken in the
usual model,
and any visit to the root in the nonbacktracking model corresponds to
a visit in the usual model.
\end{pf}

For the remainder of this section, we only consider the nonbacktracking
frog model.
We record an observation: Suppose the initial frog in the
nonbacktracking model steps from the
root~$\varnothing$ to a child~$\varnothing'$. Since frogs are stopped
at the root, no other
child of the root besides $\varnothing'$ is ever visited, and all
action occurs in the subtree
rooted at $\varnothing'$.




\subsection{Formal definition of $\mathcal{A}$} \label{secAdef}

%
%
%
%
%
%
%

\begin{figure}
\begin{tabular}{@{}cc@{}}

\includegraphics{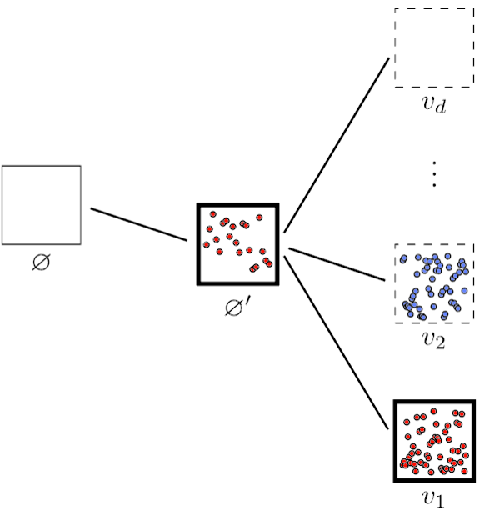}
 & \includegraphics{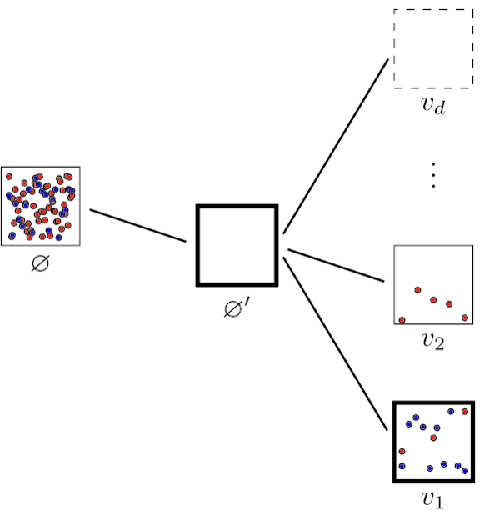}\\
\footnotesize{(a)} & \footnotesize{(b)}
\end{tabular}
\caption{An interacting particle system related to the frog model.
Initially, the number of active particles at $\varnothing'$ is
distributed as $\operatorname{Poi}(\mu)$, and the number of active
particles at
$v_1$ is distributed
according to some probability measure $\pi$. Active particles take
random nonbacktracking steps until reaching a~leaf.
For each $2\leq i\leq d$, if any of these
particles reach $v_i$, then a new $\pi$-distributed batch of particles
is released at $v_i$.
These second-wave particles do not activate other vertices.
\textup{(a)} \textit{Initial state}: particles at $\varnothing'$ and
$v_1$ will move first and possibly release a second wave of particles
from $v_2,\ldots,v_d$. \textup{(b)} \textit{Terminal state}: $\#\{
\mbox{particles at
$\varnothing$}\} \sim\mathcal{A} \pi$.}
\label{figAasystem}
\end{figure}

Fix a probability measure $\pi$ on the nonnegative integers. We will
define $\mathcal A \pi$ to be the probability measure
for the number of particles ending at $\varnothing$ (see Figure~\ref{figAasystem}) in the random system of nonbacktracking particles
described below.


The setting for the particle system is a star graph, consisting of a
central vertex
connected to $d+1$ leaf vertices. In a slight abuse of notation, we
reuse the vertex names from
Figure~\ref{figtree}, calling the central vertex $\varnothing'$ and
the leaves
$\varnothing$ and $v_1,\ldots,v_d$.
%
%
%
Let $X\sim\operatorname{Poi}(\mu)$ and $X_1,\ldots,X_d\sim\pi$,
all independent.
Place $X$ particles at $\varnothing'$ and $X_i$ particles at each $v_i$.
Each particle if activated will
perform an independent random nonbacktracking walk until it halts at a leaf.

Initially, only the particles at $\varnothing'$ and at $v_1$ are active.
If one of these first-wave particles lands at $v_i$ for $i\geq2$, then
the particles there
are activated and begin independent nonbacktracking random walks until
reaching a leaf.
These second-wave particles do not activate other particles; only the first-wave
particles have that power.
The number of particles that finish at $\varnothing$ in this system is
a random variable,
and we define $\mathcal{A}\pi$ as its law.
With these dynamics, we can summarize the system as follows:
\begin{itemize}
\item Particles at $\varnothing'$ move to one of $\{\varnothing,
v_1,\ldots, v_d\}$ each with probability $1/(d+1)$.
\item Particles at $v_1$ move to one of $\{\varnothing, v_2, \ldots,
v_d\}$ each with probability $1/d$.
\item If a first-wave particle visits $v_i$, the particles at $v_i$
move to $\varnothing$ with probability $1/d$.
\end{itemize}

For $2\leq i\leq d$, let $E_i$ be the event that a first-wave particle
ends at $v_i$. The following lemma follows from the definition of
$\mathcal A$. Informally, it says that conditional on how many of the
events $E_2,\ldots,E_d$ occur,
the number of second-wave particles ending at $\varnothing$
is a sum of independent thinned copies of $\pi$.

\begin{lemma}\label{lemfirstwavecondition}
Conditional on $\sum_{i=2}^d \mathbf{1}_{E_i}=u$,
the number of second-wave particles ending at $\varnothing$ is
distributed as the sum
of $u$ independent $\operatorname{Bin}(\pi,1/d)$-distributed random variables.
\end{lemma}

\begin{pf}
If $E_i$ occurs then by definition a $\pi$-distributed batch of
particles is released at $v_i$. With probability $1/d$ each released
particle halts at $\varnothing$. As particles move independently, the
total number is distributed as $\operatorname{Bin}(\pi, 1/d)$. Since
the second-wave
particles cannot wake other sites, the total number of particles to
arrive is distributed as claimed.
\end{pf}


Now, we show the connection between this operator and the frog model.


\begin{lemma}\label{lemfixedpoint}
Let $\nu$ be the distribution of number of returns to the root
in the nonbacktracking frog model on the $d$-ary tree
with sleeping frog distribution $\operatorname{Poi}(\mu)$. Then
$\mathcal{A}\nu\preceq\nu$.
\end{lemma}
\begin{pf}
Let $\mathbb{T}_d(x)$ denote the subtree of $\mathbb{T}_d$ rooted at
a given vertex~$x$.
Recall that no children of the root other than $\varnothing'$, the child
visited by the initial frog, are ever visited.
In light of this,
it will be helpful to think of the nonbacktracking frog model
as taking place on $\varnothing\cup\mathbb{T}_d(\varnothing')$
rather than on
all of $\mathbb{T}_d$.

We say that the frogs sleeping on some vertex $v\in\mathbb{T}_d(v_1)$
\emph
{wake within $\mathbb{T}_d(v_1)$}
if there exists a chain of vertices $x_1,\ldots,x_m=v$ all in $\mathbb{T}
_d(v_1)$ such that the initial frog starting
from the root visits
$x_1$, a frog starting at $x_1$ visits $x_2$, and so on. More simply,
a frog is woken
within $\mathbb{T}_d(v_1)$ if it would have been woken even if there
were no frogs
sleeping on any vertices outside of $\mathbb{T}_d(v_1)$.

We define some random variables counting frogs that might possibly
visit the root.
Let $X \sim\operatorname{Poi}(\mu)$ be the number of frogs sleeping
on $\varnothing
'$, which are woken
by the initial frog.
Let $X_1$ be the number of frogs waking within $\mathbb{T}_d(v_1)$
that visit
$\varnothing'$.
We claim that $X_1$ is distributed as $\nu$.
Indeed, when we consider frogs as waking only if they wake within
$\mathbb{T}_d(v_1)$ and relabel the vertices
$\{\varnothing'\}\cup\mathbb{T}_d(v_1)$ as $\{\varnothing\}\cup
\mathbb{T}
_d(\varnothing')$,
we see a process identical in law to the original nonbacktracking frog model.
Call the frogs counted by $X$ and $X_1$ the \emph{first-wave frogs}.

For each $2\leq i\leq d$,
let $E_i$ be the event that some of the frogs counted by $X$ or $X_1$
move to $v_i$.
Conditional on $E_i$, arbitrarily choose one of these frogs that
visits $v_i$ and call it $f$.
We say that the frogs at $v$ are woken within $\mathbb{T}_d(v_i)$ if there
exists a chain of
vertices $x_1,\ldots,x_m=v$ in $\mathbb{T}_d(v_i)$ such that $f$ visits
$x_1$, a frog
starting at $x_1$ visits $x_2$, and so on.
Let $X_i$ be the number
of frogs waking within $\mathbb{T}_d(v_i)$
that visit $\varnothing'$.
By the same argument showing that $X_1\sim\nu$, the distribution of
$X_i$ conditional on $E_i$ is also $\nu$.
Furthermore, for any $\{i_1,\ldots,i_k\}\subseteq\{2,\ldots,d\}$, the
random variables
$X_{i_1},\ldots,X_{i_k}$ are conditionally independent
given $E_{i_1},\ldots,E_{i_k}$, since each $X_i$ is determined
solely by the paths of the frogs sleeping in $\mathbb{T}_d(v_i)$.
We call the frogs counted by $X_2,\ldots,X_d$ the \emph{second-wave frogs}.

The first- and second-wave frogs all visit $\varnothing'$.
We define $V''$ as the number of these that move from there to
$\varnothing$.
\begin{claim*}
$V''\sim\mathcal{A}\nu$.
\end{claim*}
\begin{pf}
Our strategy is to show that the first-wave frogs
behave identically as the first-wave particles, and then to show
that the second-wave frogs conditional on the behavior of first-wave
frogs behave
the same as the second-wave particles conditional on the behavior of
the first-wave particles.

For the first of these claims, consider the first-wave frogs,
counted by $X$ and $X_1$. Observe that $X$ and $X_1$ are independent
with $X\sim\operatorname{Poi}(\mu)$
and $X_1\sim\nu$, just as in the particle system defining $\mathcal
{A}\nu$.
The frogs counted by $X$ move from $\varnothing'$ independently
to a random choice out of $\varnothing,v_1,\ldots,v_d$,
and the frogs counted by $X_1$ move from $\varnothing'$ independently
to a random choice out of $\varnothing,v_2,\ldots,v_d$, also matching
the particle system. Thus, the locations of the first-wave frogs
one step after leaving $\varnothing'$ are distributed identically to
the ending locations of the first-wave particles.

Now, condition on some arrangement of the first-wave frogs on
$\varnothing,v_1,\ldots,v_d$ one step after leaving $\varnothing$.
Suppose that $u$ out of the vertices $v_2,\ldots,v_d$ are occupied
by first-wave frogs in this arrangement. The number of second-wave frogs
visiting $\varnothing'$ conditional on this arrangement of first-wave frogs
is a sum of $u$ independent copies of $\nu$.
Each second-wave frog that visits $\varnothing'$ has an independent
$1/d$ chance of moving next to $\varnothing$. Thus,
the number of second-wave frogs that visit $\varnothing$ is the sum
of $u$ independent copies of $\operatorname{Bin}(\nu,1/d)$.
This matches the conditional distribution of second-wave particles
ending at $\varnothing$ given in Lemma~\ref{lemfirstwavecondition}.
Thus, the distribution of the number of first- and second-wave frogs visiting
$\varnothing$ is the same as the distribution of the number of first-
and second-wave particles ending at $\varnothing$, which is
by definition~$\mathcal{A}\nu$.
\end{pf}
With this claim, the proof of the lemma is almost complete:
Let $V$ be the total number of visits to $\varnothing$ in the nonbacktracking
frog model. Since $V''\leq V$ with $V''\sim\mathcal{A}\nu$ and
$V\sim\nu$,
we have shown that
$\mathcal{A}\nu\preceq\nu$.
\end{pf}

\subsection{Properties of $\mathcal{A}$} \label{secAproperties}

We first show \eqref{eqmonotonic}, monotonicity of $\mathcal{A}$
with respect
to stochastic dominance.

\begin{lemma}\label{lemmonotonicity}
If $\pi_1\preceq\pi_2$, then $\mathcal{A}\pi_1\preceq\mathcal
{A}\pi_2$.
\end{lemma}
\begin{pf}
If $\pi_1\preceq\pi_2$, then we can couple the two particle systems
defining $\mathcal{A}\pi_1$ and $\mathcal{A}\pi_2$
so that the second particle system contains all the same particles as
the first, moving identically,
as well as additional ones. Thus, at least as many particles visit
$\varnothing$ in the second
system as in the first, and $\mathcal{A}\pi_1\preceq\mathcal{A}\pi_2$.
\end{pf}

Now, we
describe the result of applying $\mathcal{A}$ to a Poisson
distribution, whose
thinning property
simplifies things.

\begin{lemma}\label{lemmixture}
The distribution $\mathcal{A}\operatorname{Poi}(\lambda)$ is a mixture
of Poisson distributions, given by
\begin{eqnarray}\label{eqmixture}
&& \mathcal{A}\operatorname{Poi}(\lambda) \sim\operatorname {Poi} \biggl(
\frac{(U+1)\lambda}{d} + \frac{\mu
}{d+1} \biggr),
\end{eqnarray}
where
\begin{eqnarray}\label{eqUdist}
&& U\sim\operatorname{Bin} \biggl(d-1, 1-\exp \biggl(-\frac{\lambda
}{d}-
\frac{\mu
}{d+1} \biggr) \biggr).
\end{eqnarray}
\end{lemma}

\begin{pf}
In the particle process defining $\mathcal{A}\operatorname
{Poi}(\lambda)$,
let $Y_{u\to v}$ be the number of particles that start at $u$ and
finish at $v$,
for $u\in\{\varnothing',v_1,\ldots, v_d\}$ and $v\in\{\varnothing
,v_1,\ldots,v_d\}$.
Each of the
$\operatorname{Poi}(\mu)$ particles starting at $\varnothing'$ moves
to a random neighbor.
By Poisson thinning,  the random variables $Y_{\varnothing'\to v}$ for
$v\in\{\varnothing,v_1,\ldots,v_d\}$
are independent and distributed as $\operatorname{Poi}(\mu/(d+1))$.
Similarly,
$Y_{v_1\to v}$ for $v\in\{\varnothing,v_2,\ldots,v_d\}$ are
independent and distributed
as $\operatorname{Poi}(\lambda/d)$. These two collections of random
variables are
also independent
of each other.

Thus, the number of first-wave particles that move to $v_i$
for each $2\leq i\leq d$ are independent and distributed as
$\operatorname{Poi}
(\lambda/d+\mu/(d+1) )$.
Let $U$ be the number of vertices out of
$\{v_2,\ldots,v_d\}$ that are visited. As each vertex has an independent
$1-\exp (-\lambda/d-\mu/(d+1) )$ chance of being visited, the
distribution of $U$ is
as given in \eqref{eqUdist}.
And since $U$ is determined by $Y_{\varnothing'\to v_i}$ and
$Y_{v_1\to v_i}$ for
$i=2,\ldots, d$, it is independent of $Y_{\varnothing'\to\varnothing}$
and $Y_{v_1\to\varnothing}$.

By Lemma~\ref{lemfirstwavecondition} and Poisson thinning,
the number of second-wave particles ending at $\varnothing$
is $\operatorname{Poi}(U\lambda/d)$.
The number of first-wave particles ending at $\varnothing$
is
$Y_{\varnothing'\to\varnothing}+Y_{v_1\to\varnothing}$,
independent of
$U$ and distributed
as $\operatorname{Poi} (\lambda/d+\mu/(d+1) )$.
Summing these together
yields \eqref{eqmixture}.
\end{pf}


We are nearly in a position to establish
that $\mathcal{A}^n\operatorname{Poi}(0)$ grows without limit as
$n\to\infty$.
First, we need two
technical lemmas on the Poisson distribution.


\begin{lemma}\label{lemconditionedPoi}
Let $\overline{Z}_{\lambda}$ be distributed as $\operatorname
{Poi}(\lambda)$
conditioned to be nonzero.
If $\lambda_1\leq\lambda_2$, then $\overline{Z}_{\lambda_1}\preceq
\overline{Z}_{\lambda_2}$.
\end{lemma}
\begin{pf}
Consider the Radon--Nikodym derivative of the law of $\overline
{Z}_{\lambda_2}$
with respect to the law of $\overline{Z}_{\lambda_1}$,
\begin{eqnarray*}
&& r(k)=\frac{\mathbf{P}[\overline{Z}_{\lambda_2}=k]}{\mathbf
{P}[\overline{Z}_{\lambda_1}=k]} = \frac{1-e^{-\lambda_1}}{1-e^{-\lambda_2}} e^{\lambda_1-\lambda_2} \biggl(
\frac{\lambda_2}{\lambda_1} \biggr)^k.
\end{eqnarray*}
The function $r(k)$ is increasing, and it is straightforward to show
that this implies that
$\overline{Z}_{\lambda_1}\preceq\overline{Z}_{\lambda_2}$ (or see
\cite{SS}, Theorem~1.C.1).
\end{pf}

\begin{lemma}\label{lemPoidivision}
Let $\overline{Z}^{(1)}_{\lambda/n}, \overline{Z}^{(2)}_{\lambda
/n},\ldots$ be independent
and distributed
as $\operatorname{Poi}(\lambda/n)$ conditioned to be nonzero.
Let $M$ be independent of these and be distributed as $\operatorname{Bin}
(n,1-e^{-\lambda/n})$, and let
\begin{eqnarray*}
&& Z = \sum_{i=1}^M \overline{Z}^{(i)}_{\lambda/n}.
\end{eqnarray*}
Then $Z$ is distributed as $\operatorname{Poi}(\lambda)$.
\end{lemma}
\begin{pf}
Decompose $\operatorname{Poi}(\lambda)$ as a sum of $n$ independent
copies of $\operatorname{Poi}
(\lambda/n)$.
Let $M$ be the number of these that are nonzero, and
condition on $M$ to get the desired representation.
\end{pf}

Finally, we prove \eqref{eqpoisson}.

\begin{prop} \label{proppoisson}
If $\mu> 2(d+1)\log d$, then there exists $\epsilon>0$ such that
\begin{eqnarray*}
&&\operatorname{Poi}(\lambda+\epsilon)\preceq\mathcal {A}\operatorname{Poi}(
\lambda)
\end{eqnarray*}
for all $\lambda\geq0$.
\end{prop}
\begin{pf}
Let $X\sim\operatorname{Poi}(\lambda+\epsilon)$ for some $\epsilon
>0$ to be chosen later, and
let $Y\sim\mathcal{A}\operatorname{Poi}(\lambda)$.
We start by decomposing $X$ into a sum of Poissons conditioned to be nonzero.
For any $a$, let
$\overline{Z}^{(1)}_a, \overline{Z}^{(2)}_a,\ldots$ be distributed
as $\operatorname{Poi}(a)$ conditioned
to be nonzero,
and let $Z_a\sim\operatorname{Poi}(a)$ (with no conditioning). Take
all these random variables
to be independent.
By Lemma~\ref{lemPoidivision}, we can write $X$ as
\begin{eqnarray}\label{eqXdef}
X &=& Z_{(\lambda+\epsilon)/d} + \sum_{i=1}^M
\overline {Z}^{(i)}_{(\lambda+\epsilon
)/d},
\end{eqnarray}
where
\begin{eqnarray*}
&& M\sim\operatorname{Bin} \biggl(d-1, 1-\exp \biggl(-\frac{\lambda
+\epsilon}{d} \biggr)
\biggr).
\end{eqnarray*}

We now turn to $Y$, which by Lemma~\ref{lemmixture} is distributed
as
\begin{eqnarray}\label{eqYdist}
&& \operatorname{Poi} \biggl(\frac{(U+1)\lambda}{d} + m \biggr),
\end{eqnarray}
where $m=\mu/(d+1)$ and $U\sim
\operatorname{Bin} (d-1, 1-\exp(-\lambda/d-m) )$.
Let $Y'\sim\operatorname{Poi} ((U+1)(\lambda+m)/d )$.
For each $u$, the distribution of $Y'$ conditional on $U=u$
is stochastically dominated by the distribution of $Y$ conditional
on $U=u$, simply because $\operatorname{Poi}(a)\preceq\operatorname
{Poi}(b)$ when $a\leq b$.
It follows that $Y'\preceq Y$. Thus, it suffices to show that
$X\preceq Y'$. Decomposing $Y'$ by Lemma~\ref{lemPoidivision} and using
the same notation as before, we can write $Y'$ as
\begin{eqnarray}\label{eqYpdef}
&& Y' = Z_{(\lambda+m)/d} + \sum_{i=1}^N
\overline{Z}^{(i)}_{(\lambda
+m)/d}
\end{eqnarray}
with
\begin{eqnarray*}
&& N \sim\operatorname{Bin} \biggl(U, 1-\exp \biggl(-\frac{\lambda
+m}{d} \biggr)
\biggr).
\end{eqnarray*}

These decompositions allow us to stochastically compare $X$ and $Y'$.
Assume that $\epsilon$ is chosen to be smaller than $m$.
We claim that
to show that $X\preceq Y'$,
it suffices to show that $M\preceq N$.
Indeed, we can then couple the random variables
on the right-hand sides of \eqref{eqXdef} and \eqref{eqYpdef} so that:
\begin{enumerate}[(3)]
\item[(1)] $M\leq N$;
\item[(2)] \vspace*{1.5pt}$Z_{(\lambda+\epsilon)/d}\leq Z_{(\lambda+m)/d}$;
\item[(3)] $\overline{Z}^{(i)}_{(\lambda+\epsilon)/d}\leq\overline
{Z}^{(i)}_{(\lambda+m)/d}$ for each
$i$.\label{itemZcond}
\end{enumerate}
Property~(2) is possible because $\operatorname
{Poi}(a)\preceq\operatorname{Poi}(b)$ if
$a\leq b$,
and (3) is possible by Lemma~\ref{lemconditionedPoi}.
Together, this yields a coupling of $X$ and $Y'$ with $X\leq Y'$.

Thus, it only remains to show that $M\preceq N$.
Recalling that $U$ is itself binomial, we have
\begin{eqnarray*}
N&\sim & \operatorname{Bin} \biggl(\operatorname{Bin} \biggl(d-1, 1-\exp \biggl(-
\frac{\lambda}{d}-m \biggr) \biggr), 1-\exp \biggl(-\frac{\lambda+m}{d} \biggr)
\biggr)
\nonumber
\\
&=& \operatorname{Bin} \biggl(d-1, \biggl(1-\exp \biggl(-\frac
{\lambda}{d}-m
\biggr) \biggr) \biggl(1-\exp \biggl(-\frac{\lambda+m}{d} \biggr) \biggr) \biggr).
\nonumber
\end{eqnarray*}
Since $M$ and $N$ are both binomial, proving $M \preceq N$ reduces to
comparing their parameters. The argument will be complete once we show
for some $\epsilon>0$ and all $\lambda>0$,
\begin{eqnarray}\label{eqNbound}
&& 1-\exp \biggl(-\frac{\lambda+\epsilon}{d} \biggr) \leq \biggl(1-\exp \biggl(-
\frac
{\lambda}{d}-m \biggr) \biggr) \biggl(1-\exp \biggl(-\frac{\lambda
+m}{d}
\biggr) \biggr).
\end{eqnarray}


Some basic calculus (see Lemma~\ref{lemcim} in the \hyperref[app]{Appendix}) establishes
that for all $d\geq2$,
\begin{eqnarray*}
&& e^{-2\log d} + e^{-2\log d/d} < 1.
\end{eqnarray*}
Since $m>2\log d$, we can choose $\epsilon>0$ such that
\begin{eqnarray*}
&& 1 >\exp \biggl(-\frac{\epsilon}{d} \biggr)\geq e^{-m}+e^{-m/d}.
\end{eqnarray*}
Multiplying both sides of this inequality by $e^{-\lambda/d}$ gives
\begin{eqnarray*}
&& \exp \biggl(-\frac{\lambda+\epsilon}{d} \biggr) \geq\exp \biggl(-\frac{\lambda}{d}-m
\biggr) + \exp \biggl(-\frac
{\lambda
+m}{d} \biggr).
\end{eqnarray*}
Thus,
\begin{eqnarray*}
1-\exp \biggl(-\frac{\lambda+\epsilon}{d} \biggr) &\leq &  1-\exp \biggl(-
\frac{\lambda}{d}-m \biggr) - \exp \biggl(-\frac
{\lambda
+m}{d} \biggr)
\\
&\leq & \biggl(1-\exp \biggl(-\frac{\lambda}{d}-m \biggr) \biggr) \biggl(1- \exp
\biggl(-\frac{\lambda+m}{d} \biggr) \biggr).
\end{eqnarray*}
Looking back at \eqref{eqNbound}, we have shown that $M\preceq N$.
\end{pf}
We have now proven \eqref{eqdominance}--\eqref{eqpoisson}, completing
the proof
of Proposition~\ref{proprecurrence}.



\subsection{Comparison to one-per-site results}
\label{subseccomparison}

In \cite{HJJ}, we proved that the frog model on a binary tree with one sleeping
frog per site is recurrent.
The proof has the same overarching idea as here: We use the
self-similarity of the tree to obtain a recursive distributional
relationship for
the number of returns to the root. We then use this relationship in a
bootstrapping
argument, assuming that the number of visits to the root is stochastically
larger than $\operatorname{Poi}(\lambda)$ and proving that it is in
fact stochastically
larger than $\operatorname{Poi}(\lambda+\epsilon)$.

The major difference between the two arguments is in the bootstrapping portion.
The approach in this paper using traditional stochastic domination
fails with the one-per-site frog model.
The problem is that the distributions given by successively applying
the analogue of the $\mathcal{A}$ operator in the one-per-site model
have finite
support, and hence are never stochastically greater than any Poisson
distribution.
Our proof in \cite{HJJ} instead uses an exotic definition of stochastic
dominance,
where $\pi_1$ is dominated by $\pi_2$ if the probability generating
function of
$\pi_1$ is greater than the probability generating function of $\pi_2$.

This generating function approach works better than the technique in
this paper
in some ways and worse in others. On one hand, it can handle both
deterministic and random
initial configurations.
On the other hand, the generating function approach seems confined to
small values of $d$.
It relies on a purely analytic argument that is elementary but difficult.
It seems impossible to apply this argument to an arbitrary choice of $d$.
Even for $d=3$, the generating functions to be analyzed become
extremely complicated.
The technical advance in this paper is the probabilistic argument we
give in Proposition~\ref{proppoisson},
which allows us to work on any $d$-ary tree.


%


\section{Transience}
\label{sectransience}


The main idea of our proof of transience is to consider a \emph{weight
function}
on the frog model.
To analyze the weight function, we bound the frog model by a branching
random walk.
The weight function is the
frog model analogue to a common martingale derived from branching
random walk
(see \cite{Biggins}).

\begin{prop}\label{proptransience}
If $\mathbf{E}\eta< \frac{(d-1)^2}{4d}$, then the frog model with an
independent copy of $\eta$ frogs per site on $\mathbb{T}_d$ is almost
surely transient.
\end{prop}

\begin{pf}
Let $F_n$ be the set of frogs awake at time~$n$. For $f\in F_n$, let
$\vert f\vert$
denote the level of $f$ on the tree (that is, its distance from the root).
We define a weight function
\begin{eqnarray*}
W_n &=&  \sum_{f\in F_n} e^{-\theta\vert f\vert},
\end{eqnarray*}
with $\theta$ to be chosen shortly.
Let
\begin{eqnarray*}
m &=& \frac{1}{d+1} e^{\theta} + \frac{d}{d+1} \mathbf{E}[\eta
+1]e^{-\theta}.
\end{eqnarray*}
Before we explain the meaning of this, we minimize $m$
by setting
$\theta= \log ((\mathbf{E}\eta+1)d )/2$, making
\[
m=\frac{2 \sqrt{ (\mathbf{E}\eta+1)d}}{d+1}< 1
\]
under our assumption that $\mathbf{E}\eta< \frac{(d-1)^2}{4d}$.

The strategy of the proof now is to show that $W_n\to0$, and hence that
the root eventually stops being visited. The term~$m$
gives an upper bound for the expected contribution to $W_{n+1}$
of a frog at time~$n$ in the following way:
Suppose that at time~$n$, some frog~$f$ is at level~$i$ of the tree
for any $i\geq1$.
With probability $1/(d+1)$, the next jump of $f$ is toward the root,
waking no frogs. With probability $d/(d+1)$, the jump is away from the root,
possibly waking up an $\eta$-distributed number of frogs. Thus, the
expected contribution
to $W_{n+1}$ from $f$ and any frogs it wakes at time~$n+1$ is at most
$e^{-\theta i}m$.
If $f$ is at the root at time~$n$, then
the expected contribution to $W_{n+1}$ from $f$ and the frogs it
wakes is at most $\mathbf{E}[\eta+1]e^{-\theta}$, which is bounded
by $m$
given our choice
of $\theta$. Therefore,
\begin{eqnarray*}
&& \mathbf{E}[W_{n+1}\mid W_n] \leq\sum
_{f\in F_n} e^{-\theta\vert
f\vert}m = mW_n.
\end{eqnarray*}
Thus, $W_n/m^n$ is a positive supermartingale. By the martingale
convergence theorem,
it converges almost surely to a finite limit.
Since $m^n\to0$, we also have $W_n\to0$ a.s., which implies that
eventually no frogs are present at the root.
\end{pf}

\begin{appendix}
\section*{Appendix}\label{app}
\begin{lemma} \label{lemcim}
$x^{-2} + x^{-2/x} < 1$ for all $x \geq2$.
\end{lemma}

\begin{pf}
Let $f(x) = x^{-2} + x^{-2/x}$. First, we show the inequality holds on
the interval $[2,8]$. Since $x^{-2}$ is decreasing,
\[
f(x) \leq\tfrac{1}4 + x^{-2/x}.
\]
It is\vspace*{1.5pt} easily checked that the maximum of $x^{-2/x}$ on $[2,8]$ occurs
at $x=8$ and is
less than $\frac{3}4$.

Next, we consider $x \geq8$.
L'H\^opital's rule implies that $\lim_{x \to\infty} f(x) =1$.
Thus, it suffices to confirm that $f(x)$ is increasing on $[8,\infty)$.
We compute
\[
f'(x)= 2 x^{-(2/x)-2} \bigl(\log x- x^{(2/x)-1}-1 \bigr).
\]
For $x\geq8$, it holds that $x^{(2/x)-1}< 1$. Hence,
\begin{eqnarray*}
&& f'(x) \geq2 x^{-(2/x)-2} (\log x- 2 ),
\end{eqnarray*}
which is positive
on $[8,\infty)$ since $\log8 > 2$.
\end{pf}
\end{appendix}

\printaddresses
\end{document}